\documentclass[12pt]{elsart}
\usepackage{amsmath,amsfonts,amssymb,amsbsy}
\usepackage{natbib}
\usepackage[pdftex]{graphicx}
\begin{document}

\begin{frontmatter}
\title{ 
Atomic force microscope based indentation stiffness tomography --- An asymptotic model
}
\author{I.I.~Argatov}
\ead{ivan.argatov@oulu.fi}
\address{Engineering Mechanics Laboratory, University of Oulu, 90014 Oulu, Finland}

\begin{abstract}
The so-called indentation stiffness tomography technique for detecting the interior mechanical properties of an elastic sample with an inhomogeneity is analyzed in the framework of the asymptotic modeling approach under the assumption of small size of the inhomogeneity. In particular, it is assumed that the inhomogeneity size and the size of contact area under the indenter are small compared with the distance between them. By the method of matched asymptotic expansions, the first-order asymptotic solution to the corresponding frictionless unilateral contact problem is obtained. The case of an elastic half-space containing a small spherical inhomogeneity has been studied in detail. 
Based on the grid indentation technique, a procedure for solving the inverse problem of extracting the inhomogeneity parameters is proposed. 
\end{abstract}

\end{frontmatter}

\setcounter{equation}{0}
\section{Introduction} 
\label{1isSectionI}

In recent years, the atomic force microscope (AFM) has become an indispensable tool for the indentation based characterization \citep{Fischer-Cripps2004} of mechanical properties of living samples at the nanometric scale \citep{KasasDietler2008,Loparic_et_al_2010,Plodinec_et_al_2012,KasasLongoDietler2013}. Based on the best-fitting analysis of the force-indentation (FI) curve, the so-called indentation stiffness tomography technique was proposed by \citet{RoduitSekatskiDietler2009} for distinguishing structures of different stiffness buried into the bulk of the sample. The developed methodology assumes that the presence of an inhomogeneity (inclusion)  changes the FI curve in a deterministic way reflecting the relative hardness/softness property of the inclusion and the depth of the inclusion. The validity of this concept was verified by finite element models and was proven useful in AFM-based indentation experiments on living cells 
\citep{RheinlaenderGeisseProkschSchaffer2011}.

Three-dimensional contact problems for the special case of an elastic half-space with inhomogeneities were studied by numerical techniques in a number of publications \citep{LerouxFulleringerNelias2010,ZhouChenKeer2011}. An extensive review of works on the elastic problems for inclusions in an elastic half-space was very recently presented by \citet{ZhouChenKeer2013}. 

In the present paper, we consider the frictionless unilateral contact problem for a homogeneous linearly elastic body with a small homogeneous inclusion (with no eigenstrains). Based on the obtained first-order asymptotic solution, we address the identification problem resulting in the identification of some parameters of the small inhomogeneity. 
The dynamic problem of stress-wave identification of material defects (e.\,g., cavities, cracks, inclusions) has been the subject of extensive research \citep{AlvesAmmari2001,GuzinaNintcheuBonnet2003,BonnetConstantinescu2005,AmmariKang2006,GuzinaChikichev2007}.
At the same time, the quasi static problem of indentation stiffness tomography has its own peculiarities, and the asymptotic model proposed herein for the extracting the inhomogeneity parameters based on the grid indentation data has not been developed elsewhere. 
The objective of this study was to develop a simple mathematical model for the indentation stiffness tomography.  

The rest of the paper is organized as follows. In Section~\ref{1isSection2}, we formulate the mathematical model of frictionless indentation. By means of the method of matched asymptotic expansions, in Section~\ref{1isSection3} an approximate analytical solution is obtained. The first-order asymptotic model for the indentation test is developed in \ref{1isSection4}. The case of a small spherical inhomogeneity in an elastic half-space is studied in detail in Section~\ref{1isSection5}, where a procedure for extracting the inhomogeneity parameters is proposed.
Finally, in Sections~\ref{1isSectionD} and \ref{1isSectionC}, respectively, we outline a discussion of the results obtained and formulate our conclusions.

\section{Mathematical model}
\label{1isSection2}

Suppose a homogeneous linearly elastic body without inhomogeneity occupies a three-dimensional domain $\Omega$. Let $\omega_\varepsilon$ be a small inhomogeneity with the center at a point ${\bf x}^0$ and the diameter proportional to a small positive parameter $\varepsilon$ such that $\omega_\varepsilon$ lies within the domain $\Omega$. Thus, the domain $\Omega_\varepsilon=\Omega\setminus\bar{\omega}_\varepsilon$ will represent the reference configuration of the homogeneous deformable body with the inhomogeneity $\omega_\varepsilon$. 

In the absence of volume forces, the displacement vectors ${\bf u}=(u_1,u_2,u_3)$ in $\Omega_\varepsilon$ and ${\bf u}^0=(u_1^0,u_2^0,u_3^0)$ in $\omega_\varepsilon$ satisfy the Lam\'e differential systems
\begin{equation}
\mu\nabla_x\cdot\nabla_x{\bf u}({\bf x})+(\lambda+\mu)\nabla_x\nabla_x\cdot{\bf u}({\bf x})=0,
\quad {\bf x}\in\Omega_\varepsilon;
\label{1is(2.1)}
\end{equation}
\begin{equation}
\mu_0\nabla_x\cdot\nabla_x{\bf u}^0({\bf x})+(\lambda_0+\mu_0)\nabla_x\nabla_x\cdot{\bf u}^0({\bf x})=0,
\quad {\bf x}\in\omega_\varepsilon.
\label{1is(2.2)}
\end{equation}

We assume that the inhomogeneity $\omega_\varepsilon$ is perfectly bonded to the surrounding medium, and the continuity and equilibrium conditions along the interface $\partial\omega_\varepsilon$ are formulated as follows:
\begin{equation}
{\bf u}({\bf x})={\bf u}^0({\bf x}),\quad 
\mbox{$\boldsymbol\sigma$}^{(n)}({\bf u};{\bf x})=\mbox{$\boldsymbol\sigma$}^{0(n)}({\bf u}^0;{\bf x}),
\quad {\bf x}\in\partial\omega_\varepsilon.
\label{1is(2.3)}
\end{equation}

The outer surface of the body $\Omega_\varepsilon$ is assumed to be decomposed into three mutually disjoint parts:
$\partial\Omega=\Gamma_u\cup\Gamma_\sigma\cup\Gamma_c$. Over $\Gamma_u$, the body is held fixed, while the surface $\Gamma_\sigma$ is assumed to be traction free, i.\,e.,
\begin{equation}
{\bf u}({\bf x})=0,\quad {\bf x}\in\Gamma_u; \quad
\mbox{$\boldsymbol\sigma$}^{(n)}({\bf u};{\bf x})=0,\quad {\bf x}\in\Gamma_\sigma.
\label{1is(2.4)}
\end{equation}

\begin{figure}[h!]
    \centering
    \includegraphics[scale=0.40]{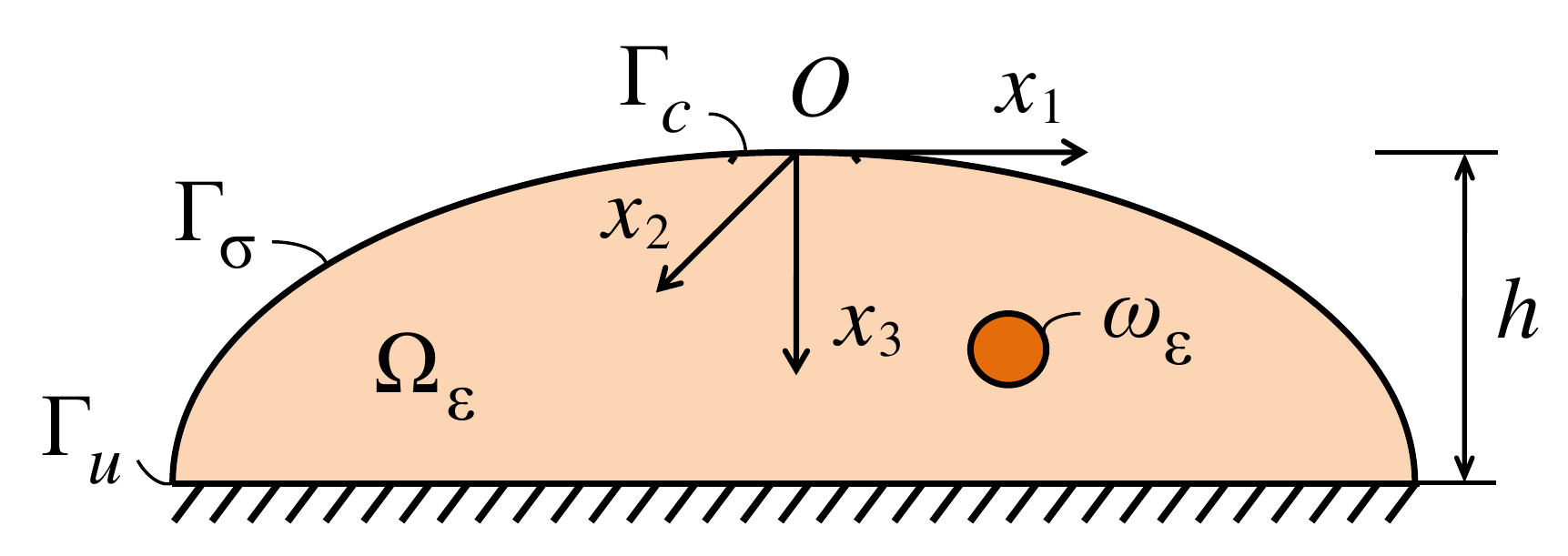}
    \caption{Deformable body and reference coordinate system.}
    \label{Fig1is1}
\end{figure}

Finally, $\Gamma_c$ is the potential contact boundary, over which contact of the deformable body $\Omega_\varepsilon$ with a rigid indenter may take place. Let us introduce a Cartesian coordinate system with the center at a point $O$ on $\Gamma_c$, which is supposed to be the initial point of contact in the case of a convex indenter (in the unloading state, the indenter touches the surface $\Gamma_c$ at the point $O$). To fix our ideas, we will assume that the $x_3$-axis is directed inside the body $\Omega_\varepsilon$ with the plane $Ox_1 x_2$ being tangent to the surface $\Gamma_c$ (see Fig.\,\ref{Fig1is1}).

\begin{figure}[h!]
    \centering
    \includegraphics[scale=0.40]{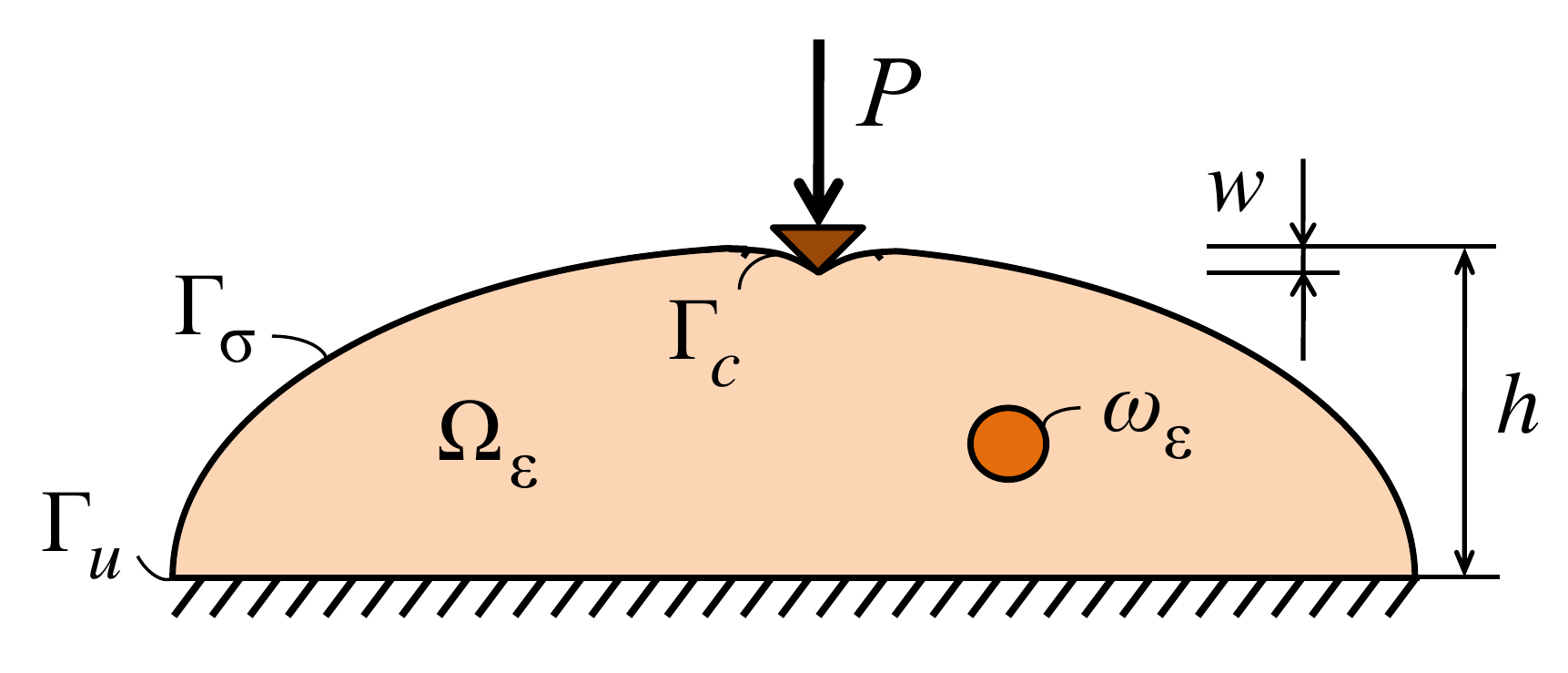}
    \caption{Indentation test schematics.}
    \label{Fig1is2}
\end{figure}

At the initial moment, the surface of the indenter is specified by the equation
\begin{equation}
x_3=-\Phi(x_1,x_2),
\label{1is(2.5)}
\end{equation}
and during the normal indentation of the indenter into the deformable body (see Fig.\,\ref{Fig1is2}), the indenter surface will be given by the equation
\begin{equation}
x_3=w-\Phi(x_1,x_2),
\label{1is(2.5a)}
\end{equation}
where $w$ is a small displacement of the indenter. 

To formulate the boundary conditions on $\Gamma_c$, following \citet{{ShillorSofoneaTelega2004}}, we denote by $g_N$ the variable gap between the indenter surface and $\Gamma_c$. Note that
$$
g_N(x_1,x_2)=g_N^0(x_1,x_2)+\frac{w}{\cos({\bf n}({\bf x}),x_3)},
$$
where ${\bf n}$ is the outer unit normal vector, $g_N^0$ is the initial gap between the body and the indenter.
Thus, the Signorini boundary conditions of frictionless contact may be stated as 
\begin{equation}
\mbox{$\boldsymbol\sigma$}^{(n)}_T({\bf u};{\bf x})=0,\quad {\bf x}\in\Gamma_c;
\label{1is(2.6)}
\end{equation}
\begin{equation}
\begin{array}{c}
\sigma^{(n)}_N({\bf u};{\bf x})\leq 0, \quad u_N({\bf x})-g_N(x_1,x_2)\leq 0, \\
\sigma^{(n)}_N({\bf u};{\bf x})(u_N({\bf x})-g_N(x_1,x_2))=0,\quad {\bf x}\in\Gamma_c.
\end{array}
\label{1is(2.7)}
\end{equation}
Here, $\sigma^{(n)}_N$ and $\mbox{$\boldsymbol\sigma$}^{(n)}_T$ are the normal and tangential components of the stress vector $\mbox{$\boldsymbol\sigma$}^{(n)}$.

Under the assumption that the size of the contact area is small compared to the curvature radii of the surface $\Gamma_c$ at the point $O$, we will have 
\begin{equation}
g_N(x_1,x_2)\simeq\Phi(x_1,x_2)-w,
\label{1is(2.8)}
\end{equation}
where $\Phi(x_1,x_2)$ is the indenter shape function introduced by Eq.~(\ref{1is(2.5)}).

Local equilibrium at the contact interface in the $x_3$-direction requires that the contact pressure should be balanced by the contact force, $P$, that is
\begin{equation}
P=-\iint\limits_{\Gamma_c}\sigma^{(n)}_N({\bf u};{\bf x})\cos({\bf n}({\bf x}),x_3)\,dS_x.
\label{1is(2.9)}
\end{equation}

Now, combining (\ref{1is(2.1)})\,--\,(\ref{1is(2.9)}), we formulate a mathematical model of unilateral contact for the linearly elastic body $\Omega_\varepsilon$ with the small inhomogeneity $\omega_\varepsilon$. 

\section{Asymptotic approximation for the displacement field away from the indenter}
\label{1isSection3}

\begin{figure}[h!]
    \centering
    \includegraphics[scale=0.40]{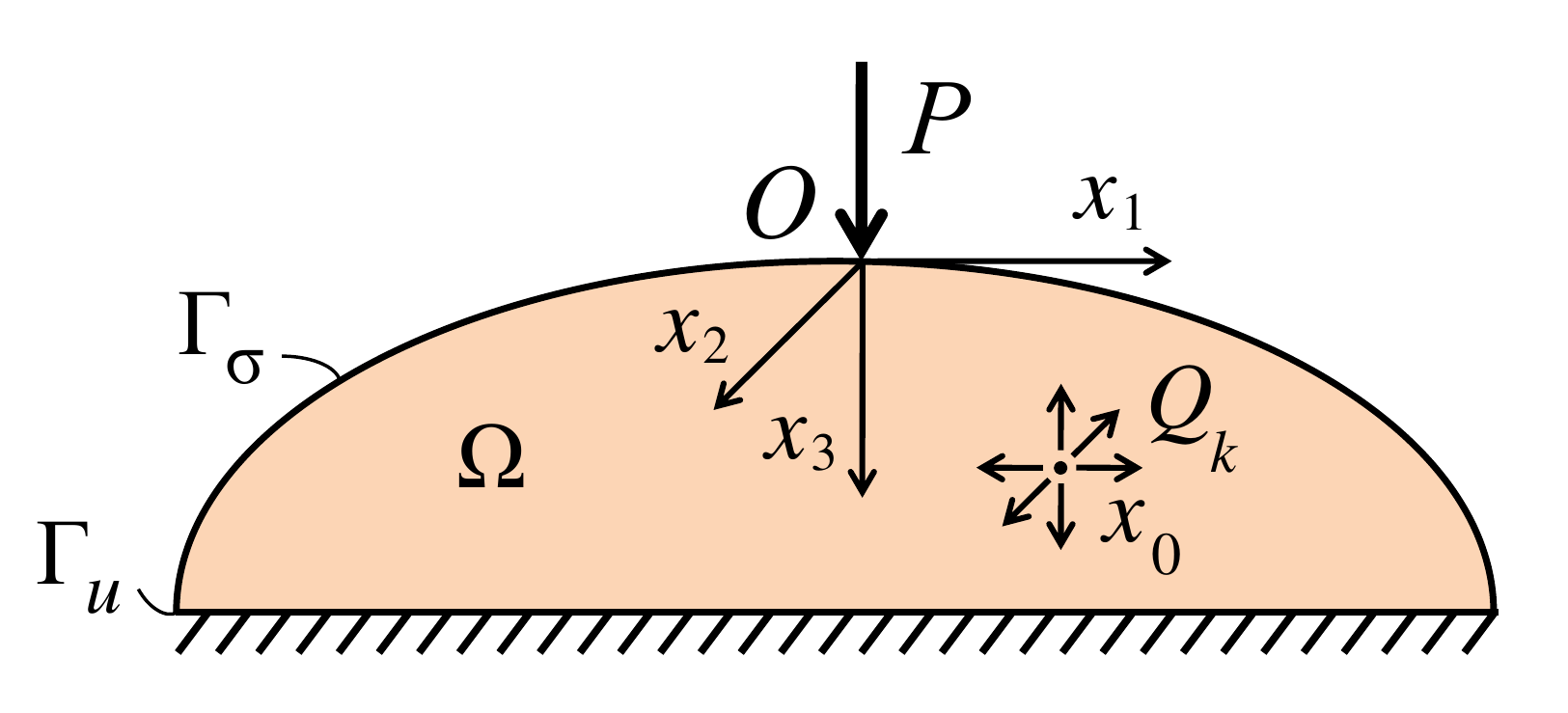}
    \caption{Schematics for the ``outer'' displacement field.}
    \label{Fig1is3}
\end{figure}

At a distance from the contact zone as well as from the inhomogeneity, the stress-strain state  of the body $\Omega_\varepsilon$ is approximated by the following so-called ``outer'' displacement vector-field (see Fig.\,\ref{Fig1is3}):
\begin{equation}
{\bf v}({\bf x})=P{\bf G}({\bf x})+\sum_{k=1}^6 Q_k {\bf M}^{(k)}({\bf x}^0,{\bf x}).
\label{1is(3.1)}
\end{equation}
Here, ${\bf G}({\bf x})$ is the solution of the elastic problem in the domain $\Omega$ of the action of a unit point force applied at the point $O$ in the $x_3$-direction, ${\bf M}^{(k)}({\bf x}^0,{\bf x})$ are the solutions of the elastic problem in the domain $\Omega$ of the action of force dipoles at the point ${\bf x}^0$. 

For the simplicity sake let us assume that the principal curvatures of the surface $\Gamma_c$ at the point $O$ are zero. Then, the vector-function ${\bf G}$  satisfies the following equations:
$$
\mu\nabla_x\cdot\nabla_x{\bf G}({\bf x})+(\lambda+\mu)\nabla_x\nabla_x\cdot{\bf G}({\bf x})=0,\quad {\bf x}\in\Omega;
$$
\begin{equation}
\mbox{$\boldsymbol\sigma$}^{(n)}({\bf G};{\bf x})=0,\quad 
{\bf x}\in\Gamma_\sigma\cup(\Gamma_c\setminus O); \qquad
{\bf G}({\bf x})=0,\quad {\bf x}\in\Gamma_u;
\label{1is(3.2)}
\end{equation}
$$
{\bf G}({\bf x})={\bf T}({\bf x})+O(1), \quad {\bf x}\to O.
$$
Here, ${\bf T}$ is the solution of the Boussinesq problem (see, e.\,g., \citet{{Johnson1985}}) of the action on the boundary of an elastic half-space $x_3>0$ of a unit force directed along the $x_3$-axis, i.\,e.,
\begin{eqnarray}
T_i({\bf x}) & = & \frac{1}{4\pi\mu}\Bigl(\frac{x_i x_3}{\vert{\bf x}\vert^3}
-\frac{\mu}{\lambda+\mu}\frac{x_i}{\vert{\bf x}\vert(\vert{\bf x}\vert+x_3)}\Bigr),\quad i=1,2,
\nonumber \\
T_3({\bf x}) & = & \frac{1}{4\pi\mu}\Bigl(\frac{x_3^2}{\vert{\bf x}\vert^3}
+\frac{\lambda+2\mu}{\lambda+\mu}\frac{1}{\vert{\bf x}\vert}\Bigr).
\nonumber
\end{eqnarray}

To describe the vector-functions ${\bf M}^{(k)}$, following \citet{ZorinMovchanNazarov1990}, we introduce vector polynomials 
\begin{equation}
\begin{array}{c}
{\bf V}^{(1)}({\bf x})=(x_1,0,0), \quad
{\bf V}^{(2)}({\bf x})=(0,x_2,0), \quad
{\bf V}^{(3)}({\bf x})=(0,0,x_3), \\
{\bf V}^{(4)}({\bf x})=\frac{1}{\sqrt{2}}(x_2,x_1,0), \quad
{\bf V}^{(5)}({\bf x})=\frac{1}{\sqrt{2}}(0,x_3,x_2), \quad
{\bf V}^{(6)}({\bf x})=\frac{1}{\sqrt{2}}(x_3,0,x_1).
\end{array}
\label{1is(3.3)}
\end{equation}

Let also $S({\bf x})$ be the Kelvin--Somigliana fundamental matrix, i.\,e.,
$$
S_{ij}({\bf x})= \frac{\lambda+\mu}{8\pi\mu(\lambda+2\mu)}\Bigl(\frac{x_i x_j}{\vert{\bf x}\vert^3}
+\frac{\lambda+3\mu}{\lambda+\mu}\frac{\delta_{ij}}{\vert{\bf x}\vert}\Bigr),\quad i,j=1,2,3. 
$$

Then, the vector-functions ${\bf M}^{(k)}({\bf x}^0,{\bf x})$ with singularities of order $O(\vert{\bf x}-{\bf x}^0\vert^{-2})$ at the point ${\bf x}^0$ are defined as follows:
$$
\mu\nabla_x\cdot\nabla_x{\bf M}^{(k)}({\bf x}^0,{\bf x})+(\lambda+\mu)\nabla_x\nabla_x\cdot
{\bf M}^{(k)}({\bf x}^0,{\bf x})=0,\quad {\bf x}\in\Omega\setminus{\bf x}^0;
$$
\begin{equation}
\mbox{$\boldsymbol\sigma$}^{(n)}({\bf M}^{(k)};{\bf x})=0,\quad 
{\bf x}\in\Gamma_\sigma\cup\Gamma_c; \qquad
{\bf M}^{(k)}({\bf x}^0,{\bf x})=0,\quad {\bf x}\in\Gamma_u;
\label{1is(3.4)}
\end{equation}
$$
{\bf M}^{(k)}({\bf x}^0,{\bf x})={\bf V}^{(k)}(\nabla_x){\bf S}({\bf x}-{\bf x}^0)+O(1), \quad {\bf x}\to{\bf x}^0.
$$

Further, to construct the boundary layer around the inhomogeneity $\omega_\varepsilon$ let us assume that 
$\omega_\varepsilon=\{{\bf x}\ \vert \ \varepsilon^{-1}({\bf x}-{\bf x}^0)\in\omega\}$,
where $\omega$ is a fixed domain, from which the small inhomogeneity $\omega_\varepsilon$ is rescaled by means of the so-called stretched coordinates $\mbox{$\boldsymbol\xi$}=\varepsilon^{-1}({\bf x}-{\bf x}^0)$.

In what follows, we make use of the Taylor expansion 
\begin{equation}
{\bf G}({\bf x})={\bf G}({\bf x}^0)+({\bf x}-{\bf x}^0)\times\mbox{$\boldsymbol\omega$}({\bf x}^0)
+\sum_{j=1}^6\epsilon_j^0{\bf V}^{(j)}({\bf x}-{\bf x}^0)+O(\vert{\bf x}-{\bf x}^0\vert^2),
\label{1is(3.5)}
\end{equation}
where $\mbox{$\boldsymbol\omega$}({\bf x})=(1/2)\nabla_x{\bf G}({\bf x})$, 
${\bf V}^{(j)}$ are the vector polynomials (\ref{1is(3.3)}), and the strain components $\epsilon_j^0$ are defined as
\begin{equation}
\begin{array}{c}
\epsilon_i^0=\varepsilon_{ii}({\bf G};{\bf x}^0), \quad i=1,2,3, \\
\epsilon_4^0=\sqrt{2}\varepsilon_{12}({\bf G};{\bf x}^0), \quad
\epsilon_5^0=\sqrt{2}\varepsilon_{23}({\bf G};{\bf x}^0), \quad
\epsilon_6^0=\sqrt{2}\varepsilon_{13}({\bf G};{\bf x}^0).
\end{array}
\label{1is(3.5a)}
\end{equation}

Thus, in view of (\ref{1is(3.1)}), $(\ref{1is(3.4)})_3$, and (\ref{1is(3.5)}), the following asymptotic expansion takes place:
\begin{eqnarray}
{\bf v}({\bf x}^0+\varepsilon\mbox{$\boldsymbol\xi$}) & = & P\Bigl\{{\bf G}({\bf x}^0)+
\varepsilon\mbox{$\boldsymbol\xi$}\times\mbox{$\boldsymbol\omega$}({\bf x}^0)
+\varepsilon\sum_{j=1}^6\epsilon_j^0{\bf V}^{(j)}(\mbox{$\boldsymbol\xi$})+
O(\varepsilon^2\vert\mbox{$\boldsymbol\xi$}\vert^2)\Bigr\}
\nonumber \\
{ } & { } & {}+\sum_{k=1}^6 Q_k\bigl\{
\varepsilon^{-2}{\bf V}^{(k)}(\nabla_\xi)S(\mbox{$\boldsymbol\xi$})+O(1)\bigr\}.
\label{1is(3.6)}
\end{eqnarray}

By the method of matched asymptotic expansions \citep{VanDyke1964,Ilin1992}, the stress-strain state of the body $\Omega_\varepsilon$ in the vicinity of the inhomogeneity is approximated by the vector-function 
\begin{equation}
{\bf w}\Bigl(\frac{{\bf x}-{\bf x}^0}{\varepsilon}\Bigr)=P\Bigl\{{\bf G}({\bf x}^0)+
({\bf x}-{\bf x}^0)\times\mbox{$\boldsymbol\omega$}({\bf x}^0)
+\sum_{j=1}^6\epsilon_j^0{\bf W}^{(j)}\Bigl(\frac{{\bf x}-{\bf x}^0}{\varepsilon}\Bigr)\Bigr\}.
\label{1is(3.7)}
\end{equation}
Here, ${\bf W}^{(j)}(\mbox{$\boldsymbol\xi$})$ are the unique solutions of the following elastic problem \citep{ZorinMovchanNazarov1990}:
$$
\mu\nabla_\xi\cdot\nabla_\xi {\bf W}^{(j)}(\mbox{$\boldsymbol\xi$})+(\lambda+\mu)\nabla_\xi\nabla_\xi\cdot
{\bf W}^{(j)}(\mbox{$\boldsymbol\xi$})=0,\quad \mbox{$\boldsymbol\xi$}\in\mathbb{R}^3\setminus\bar{\omega};
$$
$$
\mu_0\nabla_\xi\cdot\nabla_\xi {\bf W}^{0(j)}(\mbox{$\boldsymbol\xi$})+(\lambda_0+\mu_0)\nabla_\xi\nabla_\xi\cdot
{\bf W}^{0(j)}(\mbox{$\boldsymbol\xi$})=0,\quad \mbox{$\boldsymbol\xi$}\in\omega;
$$
$$
{\bf W}^{(j)}(\mbox{$\boldsymbol\xi$})={\bf W}^{0(j)}(\mbox{$\boldsymbol\xi$}),\quad
\mbox{$\boldsymbol\sigma$}^{(n)}({\bf W}^{(j)};\mbox{$\boldsymbol\xi$})=
\mbox{$\boldsymbol\sigma$}^{0(n)}({\bf W}^{0(j)};\mbox{$\boldsymbol\xi$}),\quad 
\mbox{$\boldsymbol\xi$}\in\partial\omega;
$$
$$
{\bf W}^{(j)}(\mbox{$\boldsymbol\xi$})={\bf V}^{(j)}(\mbox{$\boldsymbol\xi$})+o(1), \quad 
\vert\mbox{$\boldsymbol\xi$}\vert\to\infty. 
$$

The solutions of the above problem have the following expansions at the infinity \citep{ZorinMovchanNazarov1990}:
\begin{equation}
{\bf W}^{(j)}(\mbox{$\boldsymbol\xi$})-{\bf V}^{(j)}(\mbox{$\boldsymbol\xi$})=
\sum_{k=1}^6{\cal P}_{jk}{\bf V}^{(k)}(\nabla_\xi)S(\mbox{$\boldsymbol\xi$})+
O(\vert\mbox{$\boldsymbol\xi$}\vert^{-3}).
\label{1is(3.8)}
\end{equation}
Here, ${\cal P}_{jk}$ are the components of the so-called elasticity polarization matrix of the inhomogeneity $\omega$
\citep{ZorinMovchanNazarov1990,LewinskiSokolowski2003,AmmariKang2007}.

Thus, in view of (\ref{1is(3.6)}) and (\ref{1is(3.8)}), the asymptotic matching of the outer asymptotic representation (\ref{1is(3.1)}) and the inner asymptotic representation (\ref{1is(3.7)}) implies the equations
\begin{equation}
Q_k=\sum_{j=1}^6 P\epsilon_j^0{\cal P}_{jk}^{\varepsilon},
\label{1is(3.9)}
\end{equation}
where ${\cal P}_{jk}^{\varepsilon}=\varepsilon^3{\cal P}_{jk}$ are the polarization matrix components for the inhomogeneity $\omega_\varepsilon$.

\section{Asymptotic model for the indentation test}
\label{1isSection4}

To construct the boundary layer type solution around the indenter, we consider the following asymptotic expansion (see (\ref{1is(3.1)}) and $(\ref{1is(3.2)})_3$):
\begin{equation}
{\bf v}({\bf x})=P({\bf T}({\bf x})+{\bf g}(O))+\sum_{k=1}^6 Q_k {\bf M}^{(k)}({\bf x}^0,O)+O(\vert{\bf x}\vert).
\label{1is(3.10)}
\end{equation}
Here, ${\bf g}({\bf x})={\bf G}({\bf x})-{\bf T}({\bf x})$ is the regular part of Green's function ${\bf G}({\bf x})$. In asymptotic analysis of frictionless contact problems, an important role is played by the asymptotic constant $g_3(O)$, which can be normalized as follows \citep{Argatov2010}:
\begin{equation}
g_3(O)=-\frac{1-\nu}{2\pi\mu}\frac{a_0}{h},
\label{1is(3.11)}
\end{equation}
where $a_0$ is a dimensionless quantity, which depends on Poisson's ratio $\nu$, $h$ is a characteristic size of the domain $\Omega$ (see Fig.~\ref{Fig1is1}).

Following \citep{Argatov1999}, the first order asymptotic model of unilateral contact can be formulated in the form
\begin{equation}
p(x_1,x_2)>0,\quad (x_1,x_2)\in\Sigma_c; \qquad
p(x_1,x_2)=0,\quad (x_1,x_2)\in\partial\Sigma_c;
\label{1is(3.12a)}
\end{equation}
\begin{eqnarray}
\iint\limits_{\Sigma_c}T_3(x_1-y_1,x_2-y_2,0)p({\bf y})\,d{\bf y} & = & w-Pg_3(O)-\sum_{k=1}^6 Q_k M^{(k)}_3({\bf x}^0,O)
\nonumber\\
{} & {} & {}-\Phi(x_1,x_2),\quad (x_1,x_2)\in\Sigma_c.
\label{1is(3.12)}
\end{eqnarray}
Here, $p(x_1,x_2)$ is the contact pressure, $\Sigma_c$ is the contact domain, which should be determined as a part of solution of the integral equation (\ref{1is(3.12)}) under the conditions (\ref{1is(3.12a)}).

Let us consider the case of blunt indenter with the shape function
\begin{equation}
\Phi(x_1,x_2)=A(x_1^2+x_2^2)^{\lambda/2}.
\label{1is(3.13)}
\end{equation}
Note that the cases $\lambda=1,2,$ and $\infty$ correspond to a conical, spherical, and cylindrical indenter, respectively.

In the axisymmetric case (\ref{1is(3.13)}), using the known analytical solution \citet{Galin2008,Sneddon1965}, we can get a closed-form solution to Eq.~(\ref{1is(3.12)}). In particular, the indenter displacement $w$, the contact force $P$, and the contact radius $a$ are related by the following two equations, which are exact consequences of Eq.~(\ref{1is(3.12)}) \citep{Argatov2011}:
\begin{equation}
w-m_3^0 P=A N_1(\lambda)a^\lambda,
\label{1is(3.14)}
\end{equation}
\begin{equation}
P=\theta_1 A\frac{\pi}{2}N_2(\lambda)a^{\lambda+1}.
\label{1is(3.15)}
\end{equation}
Here we introduced the following notation (with $\Gamma(x)$ being the gamma function):
$$
\theta_1=\frac{2E}{1-\nu^2}, \quad
N_1(\lambda)=2^{\lambda-2}\lambda\frac{\Gamma(\frac{\lambda}{2})^2}{\Gamma(\lambda)},\quad
N_2(\lambda)=\frac{2^{\lambda-1}\lambda^2}{\pi(\lambda+1)}
\frac{\Gamma(\frac{\lambda}{2})^2}{\Gamma(\lambda)}
$$
\begin{equation}
m_3^0=g_3(O)+\sum_{j,k=1}^6 \epsilon_j^0{\cal P}_{jk}^\varepsilon M^{(k)}_3({\bf x}^0,O).
\label{1is(3.16)}
\end{equation}

Equations (\ref{1is(3.14)}) and (\ref{1is(3.15)}) imply that during the depth-sensing indentation, the incremental indentation stiffness can be evaluated according to the relation
\begin{equation}
\frac{dP}{dw}=\frac{\theta_1 a}{1+m_3^0\theta_1 a},
\label{1is(3.17)}
\end{equation}
which represents the first-order asymptotic model for the indentation stiffness in the axisymmetric case. We emphasize that the coefficients of the asymptotic representation (\ref{1is(3.17)}) do not depend on the dimensionless parameter $\lambda$, which describes the indenter shape \citep{ArgatovSabina2013}.

It should be noted that Eq.~(\ref{1is(3.17)}) is not convenient from the point of view of the indentation stiffness tomography method, because the contact radius $a$ is not an easily controlled parameter. Based on the asymptotic analysis given \citet{ArgatovSabina2013}, the following first-order asymptotic model for the indentation stiffness can be derived from Eqs.~(\ref{1is(3.14)}) and (\ref{1is(3.15)}):
\begin{equation}
\frac{dP}{dw}\simeq \theta_1\Bigl(\frac{w}{A N_1(\lambda)}\Bigr)^{1/\lambda}\biggl\{
1-\theta_1 m_3^0\frac{(\lambda+2)}{(\lambda+1)}\Bigl(\frac{w}{A N_1(\lambda)}\Bigr)^{1/\lambda}
\biggr\}.
\label{1is(3.20)}
\end{equation}

Thus, if the indentation testing is performed with the same indentation depth, denoting by $S_\varepsilon$ the left-hand side of Eq.~(\ref{1is(3.20)}), we can rewrite it as follows:
\begin{equation}
S_\varepsilon\simeq S_0\Bigl\{1-m_3^0\frac{(\lambda+2)}{(\lambda+1)}S_0\Bigr\}.
\label{1is(3.21)}
\end{equation}
Here, $S_0$ is the indentation stiffness of the bulk material under the specified level of indentation, i.\,e.,
\begin{equation}
S_0=\theta_1\Bigl(\frac{w}{A N_1(\lambda)}\Bigr)^{1/\lambda}.
\label{1is(3.22)}
\end{equation}

From (\ref{1is(3.21)}), it immediately follows that 
\begin{equation}
\frac{S_\varepsilon-S_0}{S_0}\simeq -m_3^0\frac{(\lambda+2)}{(\lambda+1)}S_0.
\label{1is(3.23)}
\end{equation}

Eq.~(\ref{1is(3.22)}) assumes that near the point of indentation, the elastic body $\Omega$ is approximated by an elastic half-space. This assumption is valid only for relatively small values of the contact radius. 

Finally, using the second Betti's formula
$$
\begin{array}{c}
\displaystyle
\iiint\limits_\Omega\bigl(
{\bf G}({\bf x})\cdot{\cal L}(\nabla_x){\bf M}^{(k)}({\bf x}^0,{\bf x})-
{\bf M}^{(k)}({\bf x}^0,{\bf x})\cdot{\cal L}(\nabla_x){\bf G}({\bf x})\bigr)d{\bf x}=\qquad\qquad \\
\displaystyle
\qquad\qquad
\iint\limits_{\partial\Omega}\bigl(
\mbox{$\boldsymbol\sigma$}^{(n)}({\bf G};{\bf x})\cdot {\bf M}^{(k)}({\bf x}^0,{\bf x})
-{\bf G}({\bf x}^0,{\bf x})\cdot\mbox{$\boldsymbol\sigma$}^{(n)}({\bf M}^{(k)};{\bf x})\bigr)dS_x
\end{array}
$$
and taking into account the relations $\mbox{$\boldsymbol\sigma$}^{(n)}({\bf G};{\bf x})=\delta({\bf x}){\bf e}_3$ (for ${\bf x}$ on $\Gamma_c$ near the point $O$) and ${\cal L}(\nabla_x){\bf M}^{(k)}({\bf x}^0,{\bf x})={\bf V}^{(k)}(\nabla_x)(\delta({\bf x}-{\bf x}^0){\bf e}_1+\delta({\bf x}-{\bf x}^0){\bf e}_2+\delta({\bf x}-{\bf x}^0){\bf e}_3)$ as well as the boundary conditions in (\ref{1is(3.2)}) and (\ref{1is(3.4)}), we obtain after integration by parts that
\begin{equation}
{\bf M}^{(k)}({\bf x}^0,O)=-\epsilon_k^0,\quad k=1,2,\ldots,6.
\label{1is(3.18)}
\end{equation}
Thus, formulas (\ref{1is(3.16)}) and (\ref{1is(3.18)}) imply that
\begin{equation}
m_3^0=g_3(O)-\sum_{j,k=1}^6 \epsilon_j^0\epsilon_k^0{\cal P}_{jk}^\varepsilon,
\label{1is(3.19)}
\end{equation}
where the strain components $\epsilon_j^0$ are defined in (\ref{1is(3.5a)}).

\section{Example. Spherical inhomogeneity in an elastic half-space}
\label{1isSection5}

According to \citet{ArgatovSevostianov2011}, the polarization matrix ${\cal P}^\varepsilon$ for a spherical inhomogeneity $\omega_\varepsilon$ of radius $r_\varepsilon$ is related to the stiffness contribution tensor ${\bf N}^\varepsilon$ of the inhomogeneity by the formula 
$$
{\cal P}^\varepsilon=V\left(\begin{array}{cccccc}       
N_{1111}^\varepsilon & N_{1122}^\varepsilon & N_{1122}^\varepsilon & 0 & 0 & 0 \\  
N_{1122}^\varepsilon & N_{1111}^\varepsilon & N_{1122}^\varepsilon & 0 & 0 & 0 \\  
N_{1122}^\varepsilon & N_{1122}^\varepsilon & N_{1111}^\varepsilon & 0 & 0 & 0 \\  
0 & 0 & 0 & 2N_{1212}^\varepsilon & 0 & 0 \\  
0 & 0 & 0 & 0 & 2N_{1212}^\varepsilon & 0 \\  
0 & 0 & 0 & 0 & 0 & 2N_{1212}^\varepsilon   
\end{array}
\right),
$$
where $V$ is a certain reference volume. By using a known solution \citep{KachanovHandbook}, we will have
$$
N_{ijkl}^\varepsilon=K_s^\varepsilon\delta_{ij}\delta_{kl}+G_s^\varepsilon\Bigl(J_{ijkl}-\frac{1}{3}
\delta_{ij}\delta_{kl}\Bigr)
$$
where coefficients 
$$
K_s^\varepsilon=\frac{V_\varepsilon}{V}K_1\Bigl(1+\frac{K_1(1-\kappa)}{G}\Bigr)^{-1},\quad
G_s^\varepsilon=\frac{V_\varepsilon}{V}G_1\Bigl(1+\frac{2G_1(5-2\kappa)}{15G}\Bigr)^{-1},
$$
and where $J_{ijkl}=(\delta_{ik}\delta_{lj}+\delta_{il}\delta_{kj})/2$ are components of the fourth rank unit tensor, 
$K_0=\lambda_0+2G_0/3$ and $K=\lambda+2G/3$ are the bulk moduli of the sphere and the matrix, respectively, $G_0=\mu_0$ and $G=\mu$ are the corresponding shear moduli, $\kappa=(\lambda+G)/(\lambda+2G)$, and $V_\varepsilon=(4\pi/3)r_\varepsilon^3$ is the sphere volume. 

Thus, in the case of spherical inhomogeneity, we obtain 
\begin{equation}
{\cal P}_{ij}^\varepsilon=E V_\varepsilon p_{ij}(\alpha,\nu,\nu_0),
\label{1is(4.1)}
\end{equation}
where $\alpha=E_0/E$ is the inclusion--matrix elastic moduli ratio,
\begin{equation}
\begin{array}{c}
\displaystyle
p_{ii}=k_s+\frac{2}{3}g_s\quad (i=1,2,3),\qquad
p_{ij}=k_s-\frac{g_s}{3}\quad (i\not=j,\ i,j=1,2,3),
 \\
\displaystyle
p_{ii}=g_s\quad (i=4,5,6),\qquad
p_{ij}=0\quad (i\not=j,\ i,j=4,5,6),
\end{array}
\label{1is(4.2)}
\end{equation}
\begin{equation}
k_s=\frac{(1-\nu)(\alpha(1-2\nu)+2\nu_0-1)}{(1-2\nu)^2(\alpha(1+\nu)+2(1-2\nu_0))},
\label{1is(4.3)}
\end{equation}
\begin{equation}
g_s=\frac{15(1-\nu)(1+\nu_0-\alpha(1+\nu))}{2(1+\nu)(2\alpha(1+\nu)(5\nu-4)+(1+\nu_0)(5\nu-7))}.
\label{1is(4.4)}
\end{equation}

Note that in the cases of spherical cavity and rigid sphere we, correspondingly, have
\begin{equation}
k_s^{\rm cavity}=-\frac{(1-\nu)}{2(1-2\nu)^2},\quad
g_s^{\rm cavity}=-\frac{15(1-\nu)}{2(1+\nu)(7-5\nu)},
\label{1is(4.5)}
\end{equation}
\begin{equation}
k_s^{\rm rigid}=\frac{1-\nu}{(1+\nu)(1-2\nu)},\quad
g_s^{\rm rigid}=\frac{15(1-\nu)}{4(1+\nu)(4-5\nu)}.
\label{1is(4.6)}
\end{equation}

\begin{figure}[h!]
    \centering
    \includegraphics[scale=0.40]{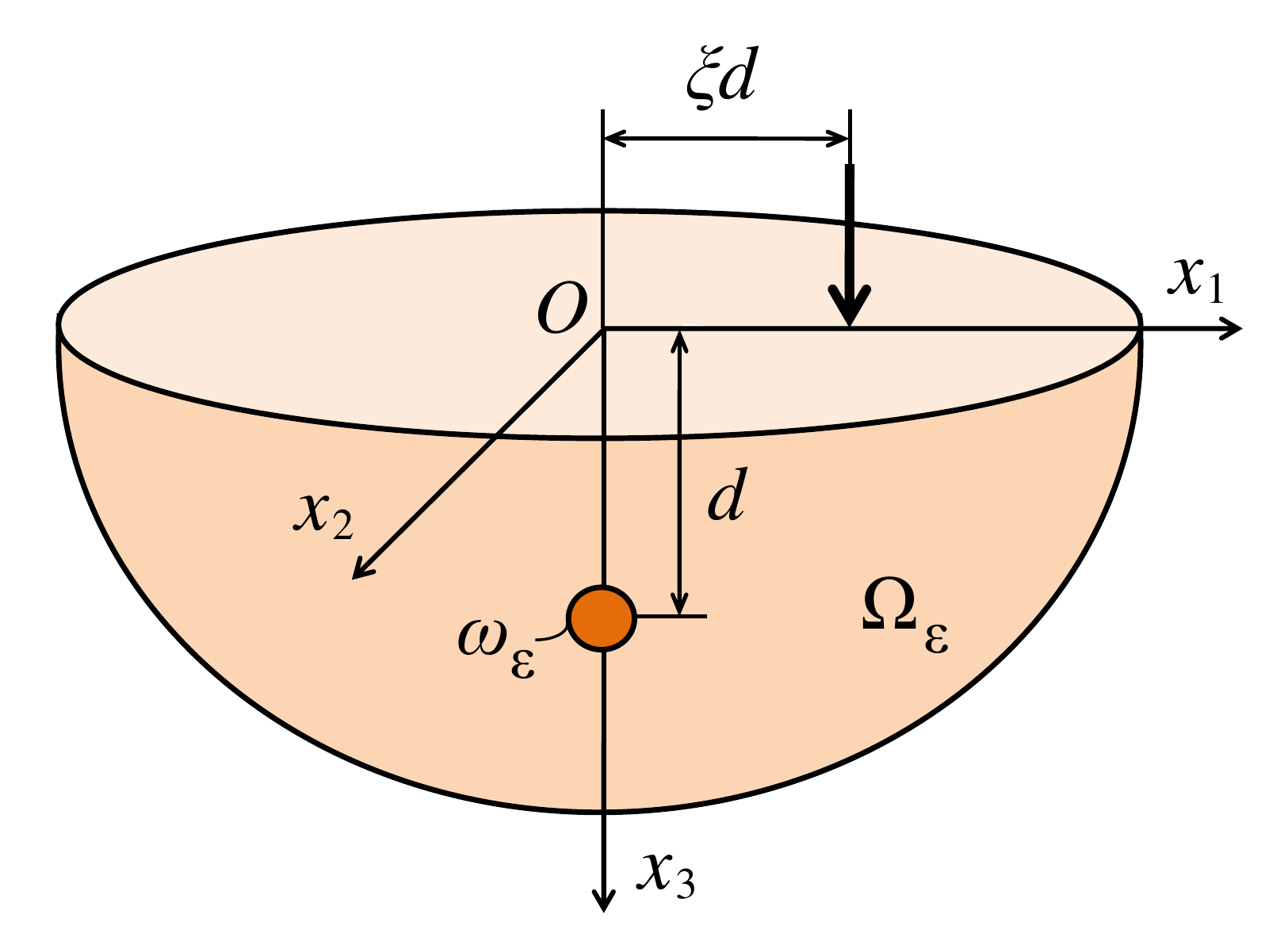}
    \caption{Schematics of the indentation test for an elastic half-space with a spherical inhomogeneity.}
    \label{Fig1is4}
\end{figure}

Now, let us introduce a coordinate system in such a way that the center of the sphere $\omega_\varepsilon$ has the coordinates $(0,0,d)$. To fix our ideas, we consider the indentation testing of the elastic half-space with the inclusion by moving the indentation point along the $x_1$-axis (see Fig.\,\ref{Fig1is4}). 

Thus, the strain components in the Boussinesq problem can be represented as
\begin{equation}
\epsilon_j^0=\frac{\bar{\epsilon}_j^0}{d^2}\quad(j=1,2,\ldots,6),
\label{1is(4.7)}
\end{equation}
where
$$
\bar{\epsilon}_1^0=\frac{1}{r^3}-\frac{3\xi^2}{r^5}+(1-2\nu)\biggl(
\frac{\xi^2}{r^2(r+1)^2}-\frac{1}{r(r+1)}+\frac{\xi^2}{r^3(r+1)}\biggr),
$$
 \begin{equation}
\bar{\epsilon}_2^0=\frac{1}{r^3}-\frac{(1-2\nu)}{r(r+1)},\quad
\bar{\epsilon}_3^0=\frac{2\nu}{r^3}-\frac{3}{r^5},
\label{1is(4.8)}
\end{equation}
$$
\bar{\epsilon}_4^0=\bar{\epsilon}_5^0=0,\quad \bar{\epsilon}_6^0=\frac{3\sqrt{2}\xi}{r^5};\quad r=\sqrt{\xi^2+1}.
$$

Because the elastic body $\Omega$ coincides with an elastic half-space, formula (\ref{1is(3.19)}) reduces to the following one:
\begin{equation}
m_3^0=-\sum_{j,k=1}^6 \epsilon_j^0\epsilon_k^0{\cal P}_{jk}^\varepsilon.
\label{1is(4.9)}
\end{equation}
Now, taking into account Eqs.~(\ref{1is(4.1)}), (\ref{1is(4.2)}), and (\ref{1is(4.7)}), we obtain
\begin{eqnarray}
m_3^0 & = & -\frac{EV_\varepsilon}{d^2}\bigl\{k_s\bigl(\bar{\epsilon}_1^0+\bar{\epsilon}_2^0+\bar{\epsilon}_3^0\bigr)^2+ g_s(\bar{\epsilon}_6^0)^2
\nonumber \\
{} & {} & {}+\frac{2}{3}g_s\bigl((\bar{\epsilon}_1^0)^2+(\bar{\epsilon}_2^0)^2+(\bar{\epsilon}_3^0)^2
-\bar{\epsilon}_1^0\bar{\epsilon}_2^0-\bar{\epsilon}_2^0\bar{\epsilon}_3^0-\bar{\epsilon}_3^0\bar{\epsilon}_1^0
\bigr)\bigr\},
\label{1is(4.10)}
\end{eqnarray}
where the expression in the braces in (\ref{1is(4.9)}) depends on the relative distance from the indentation point to the epicenter of the spherical inhomogeneity, $\xi$, the relative stiffness of the inhomogeneity, $\alpha$, and Poisson's ratios, $\nu$ and $\nu_0$.

Because a value of $0{.}5$ is usually assumed for Poisson's ratio of cells, let us consider a special case of incompressible bulk material when $\nu=0{.}5$. Then, according to Eqs.~(\ref{1is(4.8)}), we will have
\begin{equation}
\bar{\epsilon}_1^0=\frac{1}{r^3}-\frac{3\xi^2}{r^5},\quad
\bar{\epsilon}_2^0=\frac{1}{r^3},\quad
\bar{\epsilon}_3^0=\frac{2\nu}{r^3}-\frac{3}{r^5},\quad
\bar{\epsilon}_6^0=\frac{3\sqrt{2}\xi}{r^5}.
\label{1is(4.11)}
\end{equation}
Now, substituting the expressions (\ref{1is(4.3)}), (\ref{1is(4.4)}), and (\ref{1is(4.11)}) into Eq.~(\ref{1is(4.10)}) and passing to the limit as $\nu\to 0{.}5$ in the obtained result, we arrive at the formula 
\begin{equation}
-m_3^0=\frac{EV_\varepsilon}{d^2}\frac{(15\alpha-10(1+\nu_0))}{3(\alpha+1+\nu_0)}\frac{1}{(\xi^2+1)^3}.
\label{1is(4.12)}
\end{equation}

On the other hand, formula (\ref{1is(3.23)}) can be written as
\begin{equation}
\frac{(\lambda+1)}{(\lambda+2)}\frac{S_\varepsilon-S_0}{S_0^2}\simeq -m_3^0,
\label{1is(4.13)}
\end{equation}
where $S_\varepsilon$ and $S_0$ are the incremental indentation stiffnesses of the elastic half-space with the inclusion and without the inclusion under the same level of indentation.

Because the epicenter of the inclusion $\omega_\varepsilon$ is not known a priori, the so-called grid indentation technique \citep{Constantinides2006} should be employed to collect the indentation stiffness data $S_\varepsilon^i$ at points 
$(x_1^i,x_2^i)$, $i=1,2,\ldots,N$.

According to the relations (\ref{1is(4.12)}) and (\ref{1is(4.13)}), the experimental data should be fitted with the following five-parameter approximation:
\begin{equation}
{\cal S}(x_1,x_2)=S_0+\frac{C_0}{(d^2+(x_1-x_1^0)^2+(x_2-x_2^0)^2)^3}.
\label{1is(4.14)}
\end{equation}
Here, $d$ is the (unknown) depth of the inclusion, $(x_1^0,x_2^0)$ are the coordinates of its epicenter, $S_0$ is the indentation stiffness of the bulk material, and $C_0$ is some coefficient, which reflects the presence of the inclusion.

The inclusion parameter identification can be performed as follows. First, by the least square method, the five unknowns on the right-hand side of (\ref{1is(4.14)}) are determined from the solution of the optimization problem 
$$
\min_{d,x_1^0,x_2^0,S_0,C_0}\sum_{i=1}^N ({\cal S}(x_1^i,x_2^i)-S_\varepsilon^i)^2.
$$

Second, in view of (\ref{1is(4.12)})\,--\,(\ref{1is(4.14)}), the following relationship should hold true:
\begin{equation}
\frac{15\alpha-10(1+\nu_0)}{\alpha+1+\nu_0}V_\varepsilon=\frac{3(\lambda+1)}{(\lambda+2)}\frac{C_0}{E S_0^2 d^4}.
\label{1is(4.15)}
\end{equation}

Since the right-hand side of Eq.~(\ref{1is(4.15)}) is supposed to be already known, we have only one equation for determining the three parameters $\alpha$, $V_\varepsilon$, and $\nu_0$. As usual in the indentation-identification testing, we can assume some value for Poisson's ratio of the inclusion, $\nu_0$. After that, Eq.~(\ref{1is(4.15)}) allows to determine either of the two parameters $\alpha$ and $V_\varepsilon$ provided the other is known. 

\section{Discussion}
\label{1isSectionD}

First of all, observe that the elastic polarization matrix ${\cal P}^{\varepsilon}$ is positive (negative) definite for relatively hard (soft) inclusions \citep{ZorinMovchanNazarov1990}. Correspondingly, the quadratic form $\epsilon_j^0\epsilon_k^0{\cal P}_{jk}^\varepsilon$ will be sign definite. This implies that a hard inclusion produces a steeper force-displacement curve than a softer one. It should be also emphasized that though this conclusion is drawn from an asymptotic model, it has a general character \citep{RoduitSekatskiDietler2009}.

Second, the proposed identification method assumes the knowledge of both Poisson's ratios $\nu_0$ and $\nu$. Note that even in the extreme cases (cavity (\ref{1is(4.5)}) or rigid inclusion (\ref{1is(4.6)})), Poisson's ratio $\nu$ is required to extract the bulk elastic modulus $E$ from the bulk indentation stiffness $S_0$.

Third, the error of determining the elastic moduli ratio $\alpha=E_0/E$ from Eq.~(\ref{1is(4.15)}) increases when the inclusion becomes stiffer ($E_0$ increases relative to $E$).  In this case, the rigid inclusion model based on Eqs.~(\ref{1is(4.6)}) and (\ref{1is(4.10)}) should be employed to verify the obtained results. 

Finally, it is to note that the same identification procedure can be used in the case of a small inhomogeneity buried in an elastic layer provided Eqs.~(\ref{1is(4.8)}) are replaced with the corresponding equations for the strains in the elastic layer without any inhomogeneity at the point coinciding with the small inhomogeneity center.

\section{Conclusions}
\label{1isSectionC}

In this study, the problem of indentation stiffness tomography for sensing of a small inhomogeneity in a homogeneous reference elastic body is investigated via an asymptotic modeling approach. With the reference body modeled as a homogeneous elastic half-space, the inverse problem of the spherical inhomogeneity parameter identification is reduced to the least-square based minimization of a misfit between the asymptotic predictions for the indentation stiffness and the massive array of experimental data for the indentation stiffness under a specified level of indentation depth collected on the body surface by means of grid indentation technique. 

\section{Acknowledgment}
The author would like to thank Dr.~M.~Loparic (Biozentrum and the Swiss Nanoscience Institute, University of Basel) for bringing his attention to the paper of \citet{RoduitSekatskiDietler2009}.


\end{document}